\documentclass{elsart}
\usepackage{amsfonts}
\usepackage{amssymb}
\usepackage{graphicx}

\newcommand{\R}{{\Bbb R}}
\newcommand{\N}{{\Bbb N}}

\begin{document}

\begin{frontmatter}
\title{Mackey-Glass type delay differential 
equations near the boundary of  absolute stability
}
 \author[Vigo]{Eduardo Liz
},
\author[Kiev]{Elena Trofimchuk} and \author[Chile]{Sergei Trofimchuk}\\
\address[Vigo]{Departamento de Matem\'atica Aplicada II, E.T.S.I.Telecomunicaci\'on,
Universidad de Vigo, Campus Marcosende, 36200 Vigo, Spain\\ E-mail:
eliz@dma.uvigo.es}
\address[Kiev]{Department of Mathematics,
National Technical University ``KPI", Kiev, Ukraine\\ E-mail: trofimch@imath.kiev.ua}
\address[Chile]{Departamento de Matem\'aticas,
Facultad de Ciencias, Universidad de Chile, Casilla 653, Santiago, Chile\\
E-mail: trofimch@uchile.cl}

\begin{abstract} 
For an equation 
$ 
x'(t) = -x(t) + \zeta f(x(t-h)), \ x \in \R, f'(0)= -1, \ \zeta > 0,
$
with $C^3$-nonlinearity $f$ which has a negative Schwarzian 
derivative and satisfies $xf(x) < 0$ for $x\not=0$, 
we prove the convergence of all solutions to zero when 
both $\zeta -1 >0$ and $h(\zeta-1)^{1/8}$ are less than 
some constant (independent on $h,\zeta$). This result gives 
 additional insight to the conjecture about the 
equivalence between local and global asymptotical stabilities 
in the Mackey-Glass type delay differential equations.
\end{abstract}

\end{frontmatter}

%
%
%
%
%
%
%
%

\newpage

\section{Introduction and main results}
\label{sec.int}
In this note, we consider the delay differential equation
\begin{equation}
x'(t) = -x(t) + \zeta f(x(t-h)), \ x \in \R, \ \zeta > 0, 
\label{wr}
\end{equation}
where $f \in C^{3}(\R,\R)$ satisfies the following 
three basic properties \textbf{(H)}: 

{\bf (H1)} $xf(x) < 0$  for  $x\neq 0$ and $f'(0) =-1$.

{\bf (H2)} $f$ is bounded below 
and there exists at most one  point $x^*\in \R$ such that $f'(x^*)=0$. Moreover,
in this case $x^*$ is a local extremum.

{\bf (H3)} $(Sf)(x)<0$ for all $x\neq x^*$, where
$
Sf=f'''(f')^{-1}-3/2(f'')^2(f')^{-2}
$
is  the  Schwarz derivative of $f$. 

We call such a delay equation  
the Mackey-Glass type equation. 

The main purpose of this work is to give an additional insight to the following 
conjecture {\bf (C)}: {\it ``local asymptotic stability of the equilibrium $e(t)
\equiv 0$ of  Eq. (\ref{wr}) implies   global  asymptotic stability, that is,
all solutions of (\ref{wr}) converge to zero when $t$ tends to infinity".} This
conjecture was first suggested by H. Smith  (see
\cite{gyt,sm}) for Nicholson's equation, while  the above form {\bf (C)} has been
proposed in
\cite{lrt}.  Moreover, the celebrated Wright conjecture 
\cite{HMO,kuang,lprtt,lrt,ltt,wright} can be viewed as a limit case of {\bf
(C)}.  It should be noted here that  the asymptotic stability of the 
linearized equation 
\begin{equation}
x'(t) = -x(t) - \zeta x(t-h), \ x \in \R,
\label{wrll}
\end{equation}
is well studied (see \cite{hl} and Proposition \ref{hlo} below), 
while there are only few results about  the
global stability of (\ref{wr}) (e.g. see \cite{gyt,lrt} for more references). 

To formulate a criterion of asymptotical stability for Eq. (\ref{wrll}),
we define new parameters $\mu = 1/\zeta \geq 0, \ \nu = \exp(-h)/\zeta  \geq 0$.
\begin{prop}[\cite{hl}]
\label{hlo}
Suppose that $\mu \geq 1$, or $\mu < 1$ and 
\begin{equation}
\label{form1}
\nu > \nu_1(\mu) = \mu \exp (\frac{-\mu\arccos(-\mu)}{\sqrt{1-\mu^2}}).
\end{equation}
Then Eq. (\ref{wrll}) is uniformly exponentially stable. 
\end{prop}
Next, the following global stability result was proved in \cite{gyt}: 
\begin{prop}
\label{main}
Assume that $f$ satisfies hypotheses {\bf (H)}. If $\mu \geq 1$, 
or $\mu < 1$ and 
\begin{equation}
\label{form}
\nu\geq \nu_2(\mu) = \frac{\mu-\mu^2}{1+\mu^2} 
\end{equation}
then the steady state $e(t)\equiv 0$ attracts all solutions $x(t)$ 
of Eq. (\ref{wr}): $x(t) \to 0$ as $t \to +\infty$. 
\end{prop}
\begin{figure}
\centering
\includegraphics{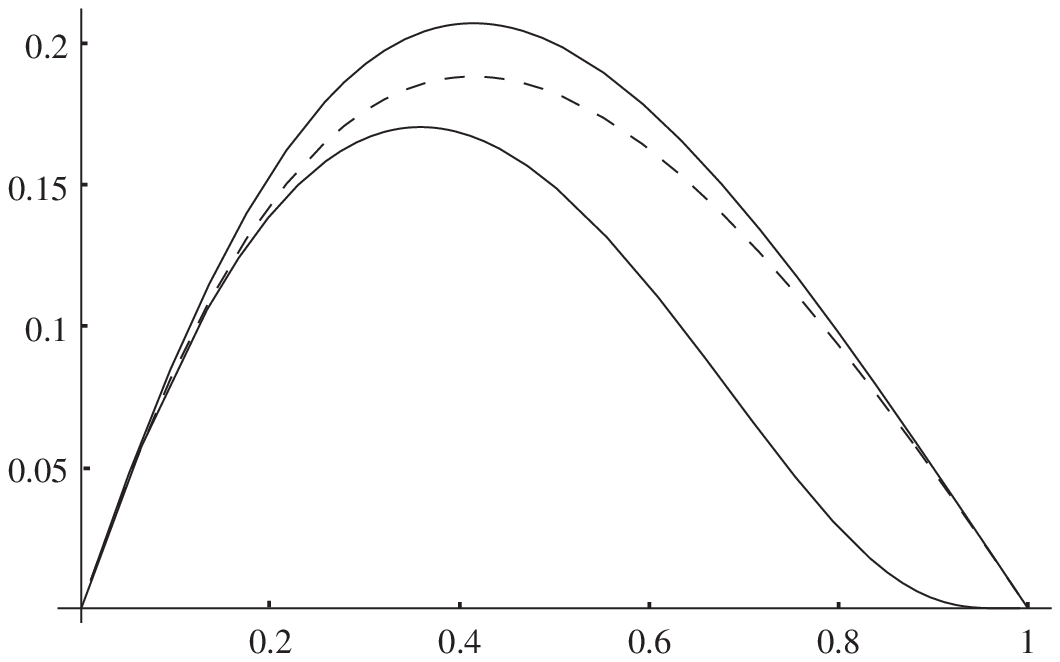}
\caption{Domains of global and local stability}
\label{fig}
\end{figure}
\begin{rem} To our best knowledge, the global stability condition 
(\ref{form}) (formulated for the Mackey-Glass type Eq.(\ref{wr})) seems to be the
best result ever reported in the literature. 

The two solid lines in Figure 1
represent the boundaries  of local and global  stability regions described in
Propositions
\ref{hlo}
 and \ref{main}: for $\mu \in (0,1)$, they are determined by the 
functions  $\nu = \nu_1(\mu)$ and $\nu = \nu_2(\mu)$ (where 
$\nu_2(\mu) > \nu_1(\mu)$). 
\end{rem}

From Fig. \ref{fig}, we observe that there is a rather
good agreement between the solid curves for sufficiently 
large $\zeta$ (e.g., for $\zeta > 5$ that corresponds to $\mu < 0.2$),
while considerable discrepancy occurs for values $\zeta$ close 
to $\zeta =1$. This difference in the behavior of these curves 
reaches its maximum at the 
point $(\mu,\nu) = (1,0)$, where the boundary of the local stability 
domain given by (\ref{form1}) (for $\mu \leq 1$) with $C^{\infty}$-smoothness 
is continued by its other part $\nu = 0$ (for $\mu \geq 1$). Indeed, 
at the same point $(\mu,\nu) = (1,0)$ the tangent line of the global 
stability curve undergoes an abrupt change. Hence, surprisingly, in order to 
 construct a counter-example to  {\bf (C)},  we should work out
parameters $\mu, \nu$  close to $(\mu,\nu) = (1,0)$. 

Moreover, there is another fact motivating the reconsideration 
of {\bf (C)}. To see this, we first state the 
following result from \cite{ilt}: 
\begin{prop}
\label{pop} 
Let $\mu>0$ and
$0<\nu<\nu_3(\mu) = \ln[(1 +\mu)/(1+\mu^2)]$. 
Then there exists a 
 periodic function $\tau: \R \to 
[0,h]$ such that the trivial solution to 
\begin{equation}
\label{nonaq}
x'(t) = -x(t) + \zeta f(x(t-\tau(t))), \ x \in \R, \ \zeta > 0,
\end{equation}
is unstable. On the other hand, if $\nu > \nu_3(\mu)$, 
then the steady state $e(t) \equiv 0$ of the equation
\begin{equation}
\label{nona}
x'(t) = -x(t) + \xi (t) f(x(t-\tau (t))), \ x \in \R,
\end{equation}
is uniformly exponentially stable for every continuous function
$\tau:\R \to [0,h]$ and
for every $\xi \in L^{\infty}(\R,\R_+)$ with 
${\rm ess}\sup_{t\in \R}\xi(t) \leq \zeta$.
\end{prop}
\begin{rem}
The graph of the function $\nu_3(\mu)$ is
depicted in Fig. \ref{fig} by a dashed line (notice that $\nu_3(\mu)=
\ln(1 + (\mu-\mu^2)/(1+\mu^2)) \geq (\mu-\mu^2)/(1+\mu^2)= \nu_2(\mu)$). 
\end{rem}
Clearly, in view of the similarity of 
 (\ref{wr}) and (\ref{nonaq}), Proposition \ref{pop} provides  another 
reason to reconsider the  global 
asymptotical stability of (\ref{wr}) for $\nu > \nu_1(\mu)$ 
(at least in the vicinity of $\mu = 1$). 

Therefore, it is important to explain 
the difference in the behavior of solid curves pictured 
in Fig. \ref{fig}. We will show below that this 
difference is only due  to the insufficiently sharp form of the
stability conditions given in Proposition \ref{main}. Indeed, let
${\mathcal D} \subset \R_+^2$ be the set of all parameters $\mu, \nu$
for which Eq. (\ref{wr}) is globally asymptotically stable, and define  $\Gamma:\R_+
\to [0,0.25]$  by 
$ \Gamma (\mu) = \inf \{\nu \geq  0: \{\mu\} \times (\nu, +\infty) \subset {\mathcal
D}\}$. The next theorem
represents the main  result of the present note, and states that  functions $\nu_1$
and $\Gamma$
  have the same slope   at $\mu =1$. 
\begin{thm}
\label{T1}There exist $\epsilon = \epsilon_f > 0, \ K = K_f > 0$ 
such that Eq. (\ref{wr}) is globally stable whenever 
$0\leq \zeta - 1 \leq \epsilon$ and 
\begin{equation}
\label{ineq}
0 \leq h < K (\zeta - 1)^{-1/8}. 
\end{equation}
As a consequence, $ \Gamma$ is differentiable at $\mu=1$, and $\Gamma'(1)=0$.
\end{thm}
\begin{rem} 
$\,$
\begin{enumerate}\item[(a)] Notice that $\Gamma(\mu) \equiv 0$
for
$\mu
\geq 1$  and $0 < \Gamma(\mu) < \nu_2(\mu)$ if $\mu \in (0,1)$. Conjecture {\bf
(C)} states that $\Gamma(\mu) = \nu_1(\mu)$; however,   we are now even unable to
prove  the continuity of $\Gamma$ over the  interval $(0,1)$, although
$\Gamma$ is lower semi-continuous thanks to the robustness of global attractivity.
\item[(b)] It should be noted that, in a small neighbourhood of $(\mu,\nu)=(1,0)$,
Eq. (\ref{wr}) can be viewed as a singularly perturbed equation \cite[Section
12.7]{hl}
$$
\varepsilon x'(t) = -x(t) + \zeta f(x(t-1)),\, \varepsilon=h^{-1}.
$$ 

It is known \cite[Theorem 7.2]{hl} that  
assumptions {\bf (H)} imply the existence of $\delta>0$ such that, for every
$(\mu,\nu)\in \{(\mu,\nu)\, :\,1-\delta<\mu<1 \, ,\, 0<\nu<  \nu_1(\mu) \}$, Eq.
(\ref{wr}) has a unique slowly oscillating periodic solution with period
$T(h,\zeta)=2h+2+O(h^{-1}+|\zeta-1|)$.

\item[(c)] It can be proved  that the set ${\mathcal D}$ is open (see
\cite{HMO,walther}). If, moreover,  one can show that ${\mathcal D}$ is closed in
the metric space 
$\{(\mu,\nu)\in (0,+\infty)^2: \nu > \nu_1(\mu) 
\ {\rm for \ } \ \mu \in (0,1]\}$,  the global stability conjecture will 
be established (compare with \cite[p. 65]{HMO}). However, we do not even 
know  if ${\mathcal D}$ is simply connected (or connected). \end{enumerate}
\end{rem}
Theorem \ref{T1} will be obtained as an easy consequence of 
several asymptotic estimations, one of which is stated below: 
\begin{thm}
\label{T2}
Let $v(t,h)$ be the fundamental solution of the linear delay
differential equation 
\begin{equation}
\label{lie}
x'(t) = - x(t) - x(t-h). 
\end{equation}
Then, for every $\alpha >2$, there exist $h_0 = h_0(\alpha) > 0, \ c  = c(\alpha)> 0$
such that 
\begin{equation}
\label{es}
|v(t,h)| \leq ch\exp(- \frac{\pi^2t}{\alpha h^3}), \ t \geq 0
\end{equation}
for all $h \geq h_0$. 
\end{thm}
\begin{rem} 
\label{ese}
$\,$
\begin{enumerate}\item[(a)] By  definition, $v(\cdot,h): [-h,+\infty) \to \R$ is
the solution of Eq. (\ref{lie}) satisfying $v(0,h) = 1$ and $v(s,h)= 0$ 
for all $s \in [-h,0)$. 
\item[(b)] It is not difficult to show (see also Remark
\ref{ark})  that the factor $h^{-3}$ from the exponent in the right-hand side of 
(\ref{es}) is the best possible (asymptotically). However, we can 
not say the same about $h$ before the exponential (for example, we do 
not know if $h$ could be replaced by $\ln h$). 
\item[(c)] We can take $c(\alpha) = b\alpha(\alpha-2)^{-1}$, where
$b >0$ does not depend on $\alpha$. 
\end{enumerate}
\end{rem}
Finally, we will also need  the following simple statement, 
which is an immediate consequence of Proposition \ref{main}
and the well-known results about period-doubling 
bifurcation for  one-dimensional dynamical systems 
defined by functions with negative Schwarzian (e.g., 
see \cite[p.92]{dev}): 
\begin{thm}
\label{T3}
There exist $\epsilon_1 = \epsilon_1(f)> 0, \ K_1 = K_1(f) > 0$ 
such that every bounded solution $x: \R \to \R$ of Eq. 
(\ref{wr}) satisfies the inequality 
\begin{equation}
\label{ie}
\sup_{t \in \R}|x(t)| \leq K_1(\zeta-1)^{1/2}
\end{equation}
for 
$0 \leq \zeta - 1 \leq \epsilon_1$.
\end{thm}
The paper is organized as follows. The proof of Theorem \ref{T2}, 
which is the most difficult ingredient of our note, 
can be found in the second section. In Section 3 we 
prove Theorem \ref{T3} and our main result (Theorem \ref{T1}), and in the
last section we discuss some other aspects of the global 
stability conjecture {\bf (C)}. 

\section{Proof of Theorem \ref{T2}}

\label{pr}
We will use the following representation of the fundamental 
solution
\begin{equation}
\label{cf}
v(t,h) = \lim_{T \to +\infty}\frac{1}{2\pi}
\int_{-T}^{T}\frac{\exp((c+is)t)}{p(c+is,h)}ds,
\end{equation}
where
$p(z)=p(z,h) = z +1 + \exp(-zh)$ is the characteristic quasipolynomial 
associated with Eq. (\ref{lie}) and 
$c > \max\{\Re \lambda: p(\lambda,h) =0\}$
(see \cite[Section 1.5]{hl}).
First we get an asymptotic estimate for $|p(z,h)|$ along the
vertical lines defined by $\lambda(s) = a + is, \ s \in \R$: 
\begin{lem} 
\label{L1}
Let $\alpha>2$ and define $\beta=(2\alpha+1)/(\alpha-2)>0$.
There exists $h_1=h_1(\alpha) >0$ such that 
\begin{equation}
\label{co}
|p(\lambda(s))| \geq \frac{\pi^2}{\beta h^2}
\end{equation}
for all $ s \in [0, 2\pi/h], \ a \in [-\pi^2/(\alpha h^3),0], \ h \geq h_1$. 
\end{lem}
\begin{pf}
We prove the lemma by contradiction. Let us suppose that there exist $h_k \to +
\infty$,
$s_k \in [0, 2\pi/h_k]$ and $a_k \in [-\pi^2/(\alpha h^3_k),0]$ such that 
\begin{equation}
\label{co1}
|p(a_k+is_k)| < \pi^2/(\beta h^2_k).
\end{equation}

Without loss of generality, 
we can assume that $s_kh_k \to \phi \in [0,2\pi]$ and 
$a_kh^3_k \to \psi \in [-\pi^2/\alpha,0]$ as $k \to \infty$. Since 
$$\lim_{k \to \infty}s_k = \lim_{k \to \infty} a_k = \lim_{k \to \infty} a_kh_k =
0,$$
we obtain from (\ref{co1}) that   $\lim_{k \to \infty} |p(a_k+is_k)| = |1
+\exp(-i\phi)| =0$. Hence $\phi = \pi$ and $\epsilon_k = s_kh_k - \pi \to 0$ when $k
\to \infty$. 

Now, it is easy to see that the inequality (\ref{co1}) implies 
$$
\frac{\pi^2}{\beta h_k} > |\pi +\epsilon_k + \exp(-a_kh_k)h_k\sin \epsilon_k|
$$
and 
$$
\frac{\pi^2}{\beta h^2_k} > |a_k +1 - \exp(-a_kh_k)\cos \epsilon_k|.
$$
The first of these inequalities is possible for all $k$ only if 
$h_k\epsilon_k \to - \pi$ as $ k \to \infty$. The second inequality 
 can be written  as 
$$
\pi^2/\beta > |a_kh_k^2 + h_k^2(1-\exp(-a_kh_k)) + h^2_k(1 -
\cos\epsilon_k)\exp(-a_kh_k)|;
$$
and takes the following limit form (when $k \to \infty$):
$$
\pi^2/\beta \geq |\psi + \pi^2/2| \geq \pi^2/2 -\pi^2/\alpha =
\frac{(\alpha-2)\pi^2}{2\alpha},
$$
a contradiction, proving Lemma \ref{L1}. 
\end{pf}
\begin{lem}
\label{cor1} For $\alpha > 2$,
there exists $h_2=h_2(\alpha)>0$ such that for every $h > h_2$, 
$ s \geq 2\pi/h$,  $a = -\pi^2/(\alpha h^3)$ we have 
\begin{equation}
\label{estr}
\max\{s -3, 0\} < \sqrt{(1+a)^2 +s^2} - \exp(-ah) \leq |p(\lambda(s))| < s +3.
\end{equation}
\end{lem}
\begin{pf}
We have, for $s >0$ and sufficiently large $h >0$, that
$$
|p(\lambda(s))| = | a + is + 1 + \exp(-ah)\exp(-ish)| \leq 1+ |a| + s + \exp(-ah) <
3 +s.
$$
On the other hand, by the triangular inequality, 
$$
|p(\lambda(s))| = |a + is + 1 + \exp(-ah)\exp(-ish)| \geq $$
$$
\geq | a + is + 1| -|\exp(-ah)\exp(-ish)| = \sqrt{(1+a)^2 +s^2} - \exp(-ah),
$$
the last part being positive for $s \geq 2\pi/h$ and  $h$ large enough. 
\end{pf}
\begin{cor}
\label{cor2} We have, for each $\alpha>2$ and all $h > h_1(\alpha)$, that 
$\sigma(h) = \max\{\Re \lambda: p(\lambda,h) =0\} < -\pi^2/(\alpha h^3)$.
\end{cor}
\begin{pf} It is well known that Eq. (\ref{lie}) is uniformly stable 
for every $h >0$ (see, e.g., \cite[p.154]{hl}), so that $\sigma (h) \leq 0$. 
It suffices now to apply Lemmas \ref{L1} and \ref{cor1} to complete the proof. 
\end{pf}
\begin{rem}
\label{ark}
In fact,  we claim that $\sigma(h) \sim -\pi^2/(2h^3)$ for $h >> 1$.

Indeed, we will establish below that the roots $\lambda(h)= a(h) \pm ib(h), \ a(h) < 0, \
b(h)> 0$   of $p(\lambda,h) =0$ have the following 
asymptotic representations for $h >> 1$: 
\begin{equation}
\label{coa}
\lambda_k(h) \sim - \pi^2(1 +2k)^2/(2h^3) \pm \pi(1+2k)i/h, \ k \in \{0,1,2,...\}.
\end{equation}
Moreover, it is easy to prove that, for $h>1$,  there exists a unique
pair of conjugate roots $\lambda (h)$ such that $|\Im (\lambda
(h))|h\leq\pi$. Thus, from (\ref{coa}) we have that, for large $h$, $\sigma(h) = \Re
(\lambda (h))\sim -
\pi^2/(2h^3)$, proving our claim.

To establish (\ref{coa}), we observe that, due to the implicit function
theorem,
$\lambda(h)$ depends  smoothly on the positive parameter $h \geq 1$.
Therefore,  rewriting the characteristic equation in the form 
\begin{eqnarray}
\label{1}
a(h) + 1 + \exp(-a(h)h)\cos(b(h)h) = 0,\\
\label{2}
b(h) = \exp(-a(h)h)\sin (b(h)h),
\end{eqnarray}
and analyzing Eq. (\ref{2}), 
we see that there exists $k \in \{0,1,2,...\}$ such that 
$b(h)h \in [2\pi k, \pi +2\pi k]$ for all $h \geq 1$ 
(notice that the characteristic equation has no 
real roots for $h \geq 1$). This means that $\lim_{h \to \infty}b(h) = 0$,
so that, by (\ref{2}), $b(h)h \to 2\pi k$
or $b(h)h \to \pi + 2\pi k$ when $h \to \infty$. We claim that 
$b(h)h \to \pi + 2\pi k$. Indeed $b(h)h \to 2\pi k$
and Eq. (\ref{1}) imply that $a(h) < -1$ for large $h$. 
This estimate allows us to conclude, again due to (\ref{1}), that 
$\lim_{h \to +\infty}a(h) = - \infty$ so that
$$\lim_{h \to +\infty}h = \lim_{h \to +\infty}|a(h)|^{-1}\ln[|a(h)+1|/\cos(b(h)h)] =0, $$
a contradiction. 
Thus $b(h)h = \pi(1+2k)+ e(h)$, where $e(h) \to 0$ as $h \to \infty$. 
Using this representation of $b(h)h$, we rewrite  Eqs. (\ref{1})
and (\ref{2}) as 
\begin{eqnarray}
\label{3}
a(h) + 1  = \exp(-a(h)h)\cos(e(h)),\\
\label{4}
\pi(1+2k)+ e(h) = - \exp(-a(h)h) h\sin (e(h)).
\end{eqnarray}
Now, Eq. (\ref{3}) implies that 
$c(h) = a(h)h \to 0$ for $h \to \infty$. Therefore, by (\ref{4}), 
we get $e(h)h = -\pi(1+2k) + o(1/h)$.  Finally, Eq. (\ref{3}) gives
$c(h)(1+o(1))= - (e^2(h)/2) (1+o(1)),$ and therefore 
$$a(h)=\frac{c(h)}{h}\sim\frac{ - e^2(h)}{2h}\sim \frac{- \pi^2(1
+2k)^2}{2h^3}.$$
\end{rem}

\begin{lem} 
\label{L2}
For each $\alpha >2$, there exist $h_0=h_0(\alpha) >0$ and $K_2=K_2(\alpha ) >0$ such
that, for every
$h > h_0$, we have
\begin{equation}
\label{coo}
\left|\lim_{T \to +\infty}\int_{-T}^{T}
\frac{e^{ist}ds}{p(-\pi^2/(\alpha h^3) +is)}\right|
\leq K_2h. 
\end{equation}
\end{lem}
\begin{pf}
First notice that the value of the integral is a real number, 
so that we have to consider only the real part of the integrand 
$e^{ist}/q(s)$. 
Since this real part is an even function, it suffices to prove 
that 
$$|I_1| =\left|\int_0^{2\pi/h}\Re [e^{ist}/q(s)]ds\right| \leq K_3 h, $$
$$|I_2| =\left|\int_{2\pi/h}^1\Re [e^{ist}/q(s)]ds\right| \leq K_4 h, 
$$ and 
$$|I_3| = \left|\int_{1}^{+\infty}\Re[e^{ist}/q(s)]ds\right| \leq K_5 h$$
for some $K_3, K_4, K_5 > 0$ and sufficiently large $h$. 

Now, by Lemma \ref{L1}, we have that, for $h\geq h_1$,
$$|I_1| \leq \int_0^{2\pi/h}|q(s)|^{-1}ds \leq (2\pi/h) (\pi^2/(\beta h^2))^{-1} =
2\beta h/\pi=K_3h,$$
where $\beta=(2\alpha +1)/(\alpha-2)$.

 Next, by Lemma \ref{cor1}, 
$$|I_2| \leq \int_{2\pi/h}^1|q(s)|^{-1}ds \leq \int_{2\pi/h}^1\frac{ds}{
|\sqrt{a^2(h)+s^2}- b(h)|},
$$
where $a(h) = 1 - \pi^2/(\alpha h^3)$ and $b(h) = \exp(\pi^2/(\alpha h^2)$. 
For sufficiently large $h$ and $s \in [2\pi/h, 1)$, 
we have 
$$\sqrt{a^2(h)+s^2}- b(h) > 0, \quad 1/a(h)<1+\pi/h, $$
$$\frac{2\pi}{a(h)h} >
\frac{\pi}{h},\quad \sqrt{1+(s+\pi/h)^2}+ b(h)/a(h) \leq 3,$$
so that
\begin{eqnarray*}|I_2| &\leq& \int_{2\pi/h}^1\frac{ds}{\sqrt{a^2(h)+s^2}- b(h)} = 
\int_{2\pi/(a(h)h)}^{1/a(h)}\frac{ds}{\sqrt{1+s^2}- b(h)/a(h)} \\
\noalign{\medskip}
&\leq& \int_{\pi/h}^{1+ \pi/h}\frac{ds}{\sqrt{1+s^2}- b(h)/a(h)} \leq 
\int_{0}^{1}\frac{ds}{\sqrt{1+(s+\pi/h)^2}- b(h)/a(h)} \\
\noalign{\medskip}
&=& \int_{0}^{1}\frac{\sqrt{1+(s+\pi/h)^2}+ b(h)/a(h)}
{1+(s+\pi/h)^2- b^2(h)/a^2(h)}ds \\
\noalign{\medskip}
 &\leq& 
\int_{0}^{1}\frac{3ds}{(s+\pi/h)^2 + (1- b^2(h)/a^2(h))}= R(h).
\end{eqnarray*}
Now, since 
$$\lim_{h\to +\infty}R(h)h^{-1} = \frac{3\sqrt{\alpha/2}}{\pi}\int_0^{\infty}
\frac{du}{(u+\sqrt{\alpha/2})^2 -1} = K_6 \in \R_+,
$$
we obtain that $|I_2| \leq (K_6+1)h = K_4h$ for sufficiently large $h$. 

Next, 
\begin{eqnarray*}
I_3 = \int_{1}^{+\infty}\frac{\cos (s(t+h))\exp(\pi^2/(\alpha h^2))
+\cos (st)(-\pi^2/(\alpha h^3) +1)}{|q(s)|^2}ds + 
\\
+\int_{1}^{+\infty}\frac{s\sin(st)}{|q(s)|^2}ds = I_4 +I_5. 
\end{eqnarray*}
Now,  for large $h$, 
$$
|I_4| \leq 3 \int_{1}^{+\infty}|q(s)|^{-2}ds \leq 3 
\int_{1}^{+\infty}(\sqrt{s^2 +0.75} - 1.25)^{-2}ds \leq K_7 \in \R,
$$
so that we only have to evaluate $I_5$. We obtain that 
$$
I_5 = \int_{1}^{+\infty}\frac{(s-|q(s)|)\sin(st)}{|q(s)|^2}ds
+ \int_{1}^{+\infty}\frac{\sin(st)}{|q(s)|}ds = I_6 +I_7,
$$
where, in virtue of (\ref{estr}),
$$
|I_6| \leq 3\int_{1}^{+\infty}|q(s)|^{-2}ds \leq  K_7.
$$
Finally, using Lemma \ref{cor1} again, we get
\begin{eqnarray*}
\left|I_7 - \int_{1}^{+\infty}\frac{\sin(st)}{s}ds\right| 
\leq \int_{1}^{+\infty}\frac{|s - |q(s)||\, |\sin(st)|}{s|q(s)|}ds\leq\\
\leq 3\int_{1}^{+\infty}|sq(s)|^{-1}ds \leq 3\int_{1}^{+\infty}\frac{1}{
s|\sqrt{a^2(h)+s^2}- b(h)|}ds\leq \\
\leq 3\int_{1}^{+\infty}\frac{ds}
{|\sqrt{a^2(h)+s^2}- b(h)|^{2}}
\leq
3\int_{1}^{+\infty}\frac{ds}{(\sqrt{s^2 +0.75} - 1.25)^{2}} \leq K_7,
\end{eqnarray*}
and since, for all $t \geq  0, \ h >0$, 
$$
\left|\int_{1}^{+\infty}\frac{\sin(st)}{s}ds\right| =
\left|\int_{t}^{+\infty}\frac{\sin(u)}{u}du\right| \leq \sup_{x\geq 0}
\left|\int_{x}^{+\infty}\frac{\sin(u)}{u}du\right| = K_8,
$$
we have the necessary estimate $|I_5| \leq K_9$. 
\end{pf}

\begin{pf} Now we can end the proof of Theorem \ref{T2}
noting that, by (\ref{cf}),
\begin{eqnarray*}
|v(t,h)| &\leq& \frac{\exp(-\pi^2t/(\alpha h^3))}{2\pi}
\left|\lim_{T \to +\infty}\int_{-T}^{T}\frac{\exp(ist)}{p(-\pi^2/(\alpha
h^3)+is,h)}ds\right| 
\\
&\leq& \frac{K_2}{2\pi}h \exp(-\pi^2t/(\alpha h^3))=
ch\exp(-\pi^2t/(\alpha h^3)). 
\end{eqnarray*}
\end{pf}

\section{Proof of Theorems \ref{T1} and \ref{T3}}

\label{sec.proof}
In order to prove Theorem \ref{T3}, we will need the following result:
\begin{prop}
\label{ag}
Assume that $f$ satisfies hypotheses {\bf (H1)}, {\bf (H2)} and set $f_{\zeta}=\zeta
f$. Then we have:
\begin{enumerate}
\item[(1)] The set $A_{\zeta}=\bigcap_{j=0}^{+\infty}f_{\zeta}^j(\R)$ is  global
attractor of the map $f_{\zeta}$; in particular, 
$A_{\zeta}=[a_{\zeta}, b_{\zeta}]$, and $f_{\zeta}(A_{\zeta})=A_{\zeta}$.
\item[(2)] Every bounded solution $x:\R\to\R$ to Eq. (\ref{wr}) satisfies
$a_{\zeta}\leq x(t)\leq b_{\zeta}$ for all $t\in\R$.
\item[(3)] If $A_1=\{0\}$, then  $\,\lim_{\zeta\to 1}a_{\zeta}=\lim_{\zeta\to
1}b_{\zeta}=0$.
\item[(4)] If $f'_{\zeta}(x)<0$ for all $x\in A_{\zeta}$, then
$f_{\zeta}(a_{\zeta})=b_{\zeta}$ and $f_{\zeta}(b_{\zeta})=a_{\zeta}$.
\end{enumerate}
\end{prop}
\begin{pf}
For (1) and (3), see \cite[Sections 2.4, 2.5]{HMO}; (4) is an immediate consequence
of (1), and, finally, (2) was established in \cite{gyt}.
\end{pf}

\begin{pf}[{\bf Theorem \ref{T3}}] 
First, note that
 hypotheses {\bf (H1)}, {\bf (H2)} imply the existence of some
$\delta$-neighbourhood $U$ of $x=0$, where 
$f$ is strictly decreasing:  $f_{\zeta}'(x) <0, \ x \in U$, $\zeta>0$. 
Next we claim that $A_1=\{0\}$. Indeed, since
$(f^2)'(0)= 1$ and $(f^2)''(0)= 0$ we obtain, in view of the negativity 
of $Sf^2$, that $(f^2)'''(0) < 0$. Therefore zero is an asymptotically 
stable point for $f^2$ (see, for example, \cite[p.25]{el}), and hence for $f$.  By  
\cite[Proposition 7]{gyt}, these facts guarantee the global  stability of the fixed
point $x = 0$ to $f$, that is, $A_1=\{0\}$.

Next, by Proposition \ref{ag} (3),  there
exists $\sigma > 0$ such that $A_{\zeta}\subset U$ for $0 < \zeta -1 < \sigma$.
Since $f$ is decreasing on $U$ we have, in view of Proposition \ref{ag} (4), that
$f_{\zeta}(a_{\zeta})=b_{\zeta}$,
$f_{\zeta}(b_{\zeta})=a_{\zeta}$.

Now, by   \cite[Corollary 12.8]{dev}, there exists a subset 
$U_{\beta} \subset U$ and a smooth   
function $\xi: U_{\beta} \to [1, +\infty), \ \xi(0) = 1, \
\xi'(0) = 0, \ \xi''(p)> 0$  for all $p \in U_{\beta}$,  
such that $f_{\xi(p)}^2(p) = p$ 
and $f_{\xi(p)}(p) \not= p$. Thus, for $\zeta \to 1^+$, we have that
$\zeta=\xi(p_1)=\xi(p_2)$ for some $p_1,p_2\in U_{\beta}$, $p_1<0<p_2$.
Moreover,the negativity of $Sf^2_{\zeta}$ 
and monotonicity of $f^2_{\zeta}$ inside
$U$
imply that $p_1$ and $p_2$ are the unique nonzero fixed points for $f^2_{\zeta}$
(in particular, $f_{\zeta}(p_1)=p_2$).  Hence, $p_1=a_{\zeta}$, $p_2=b_{\zeta}$,
and, by Proposition \ref{ag} (2), every bounded solution
$x: \R \to
\R$ to  Eq. (\ref{wr}) satisfies the inequality 
$$p_1 \leq x(t) \leq p_2, \  t \in \R.  $$
Finally, using the relations $\xi(p_i) = \zeta,\, i=1,2,$ and the properties
of $\xi$, we obtain (\ref{ie}) for  $\zeta -1 > 0$ sufficiently small. 
\end{pf}

\begin{pf}[{\bf Theorem \ref{T1}}]
Let $z: \R \to \R$ be a bounded solution to Eq.(\ref{wr}). 
Then $z(t)$ satisfies the following linear equation 
\begin{equation}
\label{xxa}
x'(t) = - x(t) - x(t-h) + a(t),
\end{equation}
where $a(t) = \zeta f(z(t-h))+ z(t-h)$. 
Take now $\epsilon_1 \in (0,1), \ K_1 > 0$ as indicated in Theorem \ref{T3}.
Then for $0 < \zeta -1 \leq \epsilon_1$, we have that  
\begin{eqnarray*}
|a(t)| &\leq& K_1(\zeta-1)^{1/2}\max_{|y| \leq K_1(\zeta-1)^{1/2}}|1+\zeta
f'(y)|\\ 
&\leq& K_1(\zeta-1)\max_{|y| \leq K_1(\zeta-1)^{1/2}}
[(\zeta-1)^{1/2}+\zeta K_1|f''(y)|]\leq \tilde{K}(\zeta-1), 
\end{eqnarray*}
where $\tilde{K} = K_1(1 + 2K_1\max_{|y| \leq K_1}
|f''(y)|)$. 

Since Eq. (\ref{xxa}) is asymptotically stable and $a(t)$ is bounded and 
continuous, it has a unique bounded solution $x(t) = z(t)$ defined for 
all $x \in \R$.  Moreover, 
$$
z(t) = \int_{-\infty}^{t}v(t-s,h)a(s)ds,
$$
so that, using Theorem \ref{T2} for an arbitrarily fixed
$\alpha>2$, we get 
\begin{eqnarray*}
|z(t)|\leq \tilde{K}(\zeta-1)\int_{-\infty}^{t}|v(t-s,h)|ds 
\leq \tilde{K}(\zeta-1)\int_{-\infty}^{t}ch\exp(- \frac{\pi^2(t-s)}{\alpha h^3})ds \\
= \frac{\tilde{K}c\alpha}{\pi^2}(\zeta-1)h^4 < (1/2)K_1 (\zeta-1)^{1/2},
\end{eqnarray*}
for
$h\geq h_0$ whenever $ h(\zeta-1)^{1/8} < K = \sqrt{\pi}
(K_1/(2\tilde{K}c\alpha))^{1/4}$. By repeating the same
argument, we can prove that $|z(t)| <
(1/2)^nK_1(\zeta-1)^{1/2}$ for all $t\in\R$ and
$n\in\N$. Thus
$z(t)
\equiv 0$. 

Without loss of generality, we can assume that
$h_0<K(\zeta-1)^{-1/8}$ for all $\zeta\in [1,1+\epsilon_1]$.
Hence we have shown above that Eq.(\ref{wr}) is globally stable
for
$h_0\leq h<K(\zeta-1)^{-1/8}$ and $0\leq \zeta -1\leq
\epsilon_1$. Finally, Proposition
\ref{main} permits us to find $\epsilon_2>0$ such that $0\leq \zeta -1\leq
\epsilon_2$ implies that (\ref{wr}) is globally stable for $0\leq h<h_0$. Thus
inequality (\ref{ineq}) is proved choosing $\epsilon=\min\{\epsilon_1,\epsilon_2\}$.

Now, (\ref{ineq}) implies that, for $\delta>0$ sufficiently small,
$
0\leq\Gamma(1-\delta)\leq
F(1-\delta),
$
where  $F(\mu)= \mu \exp\{-K
((1/\mu)-1)^{-1/8}\}$.

Since $\displaystyle\lim_{\delta\to
0^+}F(1-\delta)/\delta=F'(1^-)=0$, we can conclude that
$\Gamma'(1)=0.$
\end{pf}

\section{Remarks and discussion}

It is easy to see that the study of  global asymptotical 
stability of the unique positive equilibrium to 
the following well-known (e.g., see \cite{cook,gl,gyt,kuang,sm}) 
delay equations (with positive $\zeta, a, x$) 
\begin{eqnarray}
\label{mg}
x'(t) = - x(t) + \frac{\zeta a^n}{a^n + x^n(t-h)},\  
n > 1, \ {\rm (Mackey-Glass)} \\
\label{lw}
x'(t) = -  x(t) + \zeta \exp(-a x(t-h))\quad {\rm (Lasota-Wazewska)}
\end{eqnarray}
\noindent can be reduced, via a simple change (a translation) of variable,
to the investigation of  global attractivity of the trivial 
solution to Eq. (\ref{wr}). In some cases (e.g. when $\zeta$ 
is close to 1), the same observation is also valid for the 
equations 
\begin{eqnarray}
\label{mg1}
x'(t) = - x(t) + \frac{ \zeta a^n x(t-h)}{a^n + x^n(t-h)},\  
n > 1, \ {\rm (Mackey-Glass)} \\
\label{nic}
x'(t) = - x(t) + \zeta x(t-h)\exp(-a x(t-h))\quad {\rm (Nicholson)}.
\end{eqnarray}
As  mentioned before (see Propositions 
\ref{hlo}, \ref{main} and \cite{gyt}), in the domain $(h,\zeta) \in \R_+^2$, 
the decay dominance condition $1 \geq \zeta$ (or $\mu \geq 1$) 
determines all parameters for which Eq. (\ref{wr}) (and, in particular,
(\ref{mg})-(\ref{nic})) is absolutely stable (we recall here that ``absolute
stability" means  ``delay independent stability"). 

Now let $\zeta > 1$ and denote by $h_c(\zeta)$ the global
stability delay threshold: $h_c(\zeta)$ is the maximal 
positive number for which the inequality $ h < h_c(\zeta)$ 
implies  convergence of all solutions to the equilibrium. 
By the above comments, it is natural to expect that 
$h_c(\zeta) \to +\infty$ as $\zeta \to 1+$ (while 
the folklore statement: {\it ``Small delays are harmless"}
implies that $h_c > 0$). Indeed, by Proposition \ref{main},
$h_c(\zeta) \geq \ln(\zeta+\zeta^{-1}) - \ln(\zeta-1) \sim - \ln(\zeta-1)$, 
so that for every $h > 0$ the Mackey-Glass delay differential 
equation can be stabilized by choosing $\zeta > 1$ sufficiently 
close to 1. This means that even  large delays are harmless 
near the boundary of absolute stability. Moreover, 
Theorem \ref{T1} has improved the above logarithmic 
estimation of $h_c(\zeta)$ near $\zeta =1$ saying 
that we have there, for some $K >0, \epsilon > 0$,  
\begin{equation}
\label{hc}
h_c(\zeta) \geq K(\zeta-1)^{-1/8} \quad  {\rm if} \ \ 0< \zeta -1 < \epsilon.
\end{equation}

Now, let us indicate briefly some aspects of 
the considered problem which could be studied 
in the future. 

First, it seems that 
the exponent $-1/8$ in (\ref{hc}) can be
significantly improved (up to $-1/2$, if 
the global stability conjecture were true).  
Unfortunately, our method (when we 
establish some estimates for the global 
attractivity domain (Theorem \ref{T3}) and 
then use the contractivity argument inside 
this domain (Theorem \ref{T2})) does not allow 
this improvement at all. The best estimate 
within our approach is $-1/6$, and it 
could be reached if we were able to show 
that $h$ before the exponential in (\ref{es}) 
is not necessary (or at least could be 
replaced with $\ln h$, see also Remark \ref{ese}). 

Second, it will be very interesting to obtain 
some $K, \epsilon$ in (\ref{hc}) explicitly. 
Also, in the statement of Theorem \ref{T1}, 
both  constants depend on the nonlinearity 
$f$: we hope that this dependence can be 
discarded with a different  approach. 

Finally, we note that  at the moment of the acceptance of this paper  we
already proved
that the inequality $\nu > \nu_3(\mu), \ \mu \in(0,1)$
(see Proposition \ref{pop}) is also sufficient for the
global stability in (\ref{wr}) (see \cite{ltt}).

%

\section*{Acknowledgments}

 This research has been supported in part
by MCT (Spain), project BFM2001-3884, and grant
FONDECYT 8990013 (Chile). The authors wrote the
paper while E. Trofimchuk and S. Trofimchuk were visiting the University of Vigo
(Spain). It is their pleasure to thank the University for its kind hospitality.

This paper was considerably improved by the
suggestions of  an anonymous referee, and we
are deeply  indebted to him.

\end{document}